\documentclass[a4paper,3p,times]{elsarticle}
\usepackage{amsthm,amsmath,amssymb,amsfonts,subcaption,listings,graphicx,moreverb,multicol}
\usepackage{lineno,hyperref,verbatim,latexsym,diagbox,algorithm,color,parskip,filecontents}
\usepackage[all]{xy}
\definecolor{codegreen}{rgb}{0,0.6,0}
\definecolor{codegray}{rgb}{0.5,0.5,0.5}
\definecolor{codepurple}{rgb}{0.58,0,0.82}
\definecolor{backcolour}{rgb}{0.95,0.95,0.92}
\lstdefinestyle{mystyle}{
    backgroundcolor=\color{backcolour},   
    commentstyle=\color{codegreen},
    keywordstyle=\color{magenta},
    numberstyle=\tiny\color{codegray},
    stringstyle=\color{codepurple},
    basicstyle=\ttfamily\footnotesize,
    breakatwhitespace=false,         
    breaklines=true,                 
    captionpos=b,                    
    keepspaces=true,                 
    numbers=left,                    
    numbersep=5pt,                  
    showspaces=false,                
    showstringspaces=false,
    showtabs=false,                  
    tabsize=2
}
\usepackage[below]{placeins} 
\usepackage[english]{babel}
\usepackage[utf8x]{inputenc}

\newtheorem{thm}{Theorem}[section]
\newtheorem{pro}[thm]{Proposition}
\newtheorem{lem}[thm]{Lemma}
\newtheorem{rem}[thm]{Remark}
\newtheorem{defn}[thm]{Definition}

\newtheorem{cor}[thm]{Corollary}

\journal{Quaestiones Mathematicae}

\begin{document}

\begin{frontmatter}
\author[label1,label2]{A. S. Bamunoba}
\ead{alex.bamunoba@mak.ac.ug}
\author[label1,label2]{I. Ndikubwayo\corref{cor1}}
\ead{innocent.ndikubwayo@mak.ac.ug}
\address[label1]{Makerere University, Department of Mathematics, Kampala, 256, Uganda}
\fntext[label2]{This research was funded by the NORHED-II project ``Mathematics for Sustainable Development (MATH4SD), 2021-2026" at Makerere University in collaboration with the University of Dar es Salaam and University of Bergen in Norway.}
\cortext[cor1]{Corresponding author}
\title{On the location of ratios of zeros of special trinomials}

\begin{abstract}
Given coprime integers $k, \ell$ with $k > \ell \geqslant 1$ and arbitrary complex polynomials $A(z), B(z)$ with $\deg(A(z)B(z))\geqslant 1$, we consider the polynomial sequence $\{P_n(z)\}$ satisfying a three-term recurrence $P_n(z)+B(z)P_{n-\ell}(z)+A(z)P_{n-k}(z)=0$ subject to the initial conditions $P_0(z)=1$, $P_{-1}(z)=\cdots=P_{1-k}(z)=0$ and fully characterize the real algebraic curve $\Gamma$ on which the zeros of the polynomials in $\{P_n(z)\}$ lie. In addition, we show that, for any (randomly chosen) $n\in \mathbb{Z}_{\geqslant 1}$ and zero $z_0$ of $P_n(z)$ with $A(z_0)\neq 0$, at-least two of the distinct zeros of the trinomial $D(t;z_0):={A(z_0)t^{k}+ B(z_0)t^{\ell}+1} $ have a ratio that lies on the real line and / or on the unit circle centred at the origin. This reveals a previously unknown geometric property exhibited by the zeros of trinomials of the form $t^k+at^{\ell}+1$ where $a\in \mathbb{C}-\{0\}$ is such that $a^k\in \mathbb{R}$.
\end{abstract}

\end{frontmatter}


\section{Basic notions and main result} \label{sec:int}
Recursively defined polynomials have always been a subject of interest since they can be used to explain phenomena frequently occurring in mathematics, statistics, physics and engineering. For example, (recursive) polynomial sequences arise in physics and approximation theory as the solutions of certain ordinary differential equations. In particular, the Hermite polynomials are used in quantum and statistical mechanics to describe solutions of the Schr\"odinger equation for a harmonic oscillator while Laguerre polynomials are used to describe the eigenfunctions for the Schr\"odinger operator associated with the hydrogen atom, \cite{costabile}. In combinatorics, these recursive polynomials can be shown to represent certain combinatorial objects like graphs and therefore can be used to come up with new identities and generating functions, the validity of which can be proved using the combinatorial interpretation. 


Some of the above mentioned polynomial sequences are examples of sequences of orthogonal polynomials satisfying three-term recursive formulae. These three-term recursion formulae are famous since they provide a necessary condition for polynomial sequences to be orthogonal, a highly sought for property in functional representations and numerical computations. Typically, orthogonal polynomials are used as basis functions in which to expand other more complicated functions to be used in either interpolation, approximation or numerical quadrature. In addition, interesting properties such as real-rootedness of the generated polynomials can be deduced from the orthogonality property, for details \cite{GVM}. 
Now, it is almost a tradition in mathematics that, whenever one encounters a family of polynomials (preferably generated by recursions), it is of interest to study properties such as orthogonality, unimodularity, real-rootedness, stability or interlacing roots among others, see \cite{petterbranden} for further reading.

In recent years however, several authors for example \cite{BO, ndikus, KD, Tr1, Tr2, TrZu1} have interested themselves in the problem of the location of zeros of polynomials in polynomial sequences, especially those that are generated by linear recurrences. This is partly because, knowing the zeros of some polynomials in the polynomial sequence may give information concerning the zeros of other polynomials in the same sequence. In all of these studies, the authors either search for a criterion for real-rootedness of the polynomial sequence (for applications in numerical analysis and combinatorics) or explicit determination of the limiting curve on which the zeros of the polynomials lie, see for example \cite{BO, ndikus,Tr1,Tr2}. 

In this setting, one considers arbitrary non-zero complex polynomials $A(z), B(z)$ with $\deg(A(z)B(z))\geqslant 1$, two coprime integers $k, \ell$ with $k > \ell \geqslant 1$ and the polynomial sequence $\{P_n(z)\}$ satisfying a three-term recurrence 
\begin{align}\label{uutiixxxx}
P_n(z)+B(z)P_{n-\ell}(z)+A(z)P_{n-k}(z)=0,
\end{align} 
subject to the initial conditions $P_0(z)=1$, $P_{-1}(z)=\cdots=P_{1-k}(z)=0$. If we take a solution of the form $P_n(z)=(t(z))^n$ where $n\in \mathbb{Z}_{\geqslant 0}$ and $t(z)$ is a nonzero complex rational function, then substituting it in the recurrence \eqref{uutiixxxx} yields 
\begin{align*}0=(t(z))^n=B(z)(t(z))^{n-\ell}+A(z)(t(z))^{n-k}=((t(z))^{k}+B(z)(t(z))^{k-\ell}+A(z))(t(z))^{n-k}.
\end{align*}
 Since $t(z)\not\equiv 0$, it follows that, $(t(z))^{k}+B(z)(t(z))^{k-\ell}+A(z)=0$ except probably for finitely many $z\in\mathbb{C}$. The equation $t^{k}+B(z)t^{k-\ell}+A(z)=0$ (with the variable $z$ of $t$ dropped) is called the characteristic equation of recurrence in \eqref{uutiixxxx}, while the polynomial $\displaystyle \Delta(t;z):=t^{k}+B(z)t^{k-\ell}+A(z)$ is called the characteristic polynomial of recurrence in \eqref{uutiixxxx}. In \cite[Lemma 1]{ndikus}, the ordinary generating function of the polynomial sequence $\{P_n(z)\}$ generated by \eqref{uutiixxxx} is shown to be
\begin{align*}
G(t; z):=\sum_{n=0}^{\infty}P_n(z)t^n=\frac{1}{1+B(z)t^{\ell}+A(z)t^{k}}.
\end{align*} 
We shall denote the denominator of $G(t;z)$ by $D(t;z):=A(z)t^{k}+B(z)t^{\ell}+1$, the reciprocal polynomial of $\Delta(t;z)$.
\begin{rem}
We consider $k$ and $\ell$ to be coprime since the case for non-coprime integers, does not provide any new information about the distribution of the zeros of non-zero polynomials in its sequence. To see this, we take $r, s \in \mathbb{Z}_{\geqslant 1}$ with $r=\ell d$, $s=k d$, gcd$(k,\ell)=1$ and $k > \ell \geqslant 1$. The ordinary generating function $G_R(t;z)$ for the polynomial sequence $\{R_m(z)\}$ generated by the three-term recurrence $\displaystyle R_n(z)+B(z)R_{n-r}(z)+A(z)R_{n-s}(z)=0$, subject to the initial conditions $R_0(z)=1$, $R_{-1}(z)=\cdots=R_{1-s}(z)=0$ is related to $G(t;z)$ as follows: \begin{align*}
G_R(t;z)=\sum_{m=0}^\infty R_m(z)t^m =\frac{1}{1+B(z)t^{r}+A(z)t^{s}}= \frac{1}{1+B(z)w^{\ell d }+A(z)w^{kd}} =\sum_{n=0}^\infty P_n(z)t^{nd}=  G(t^d;z).
\end{align*} 
Therefore, $R_m(z)=P_n(z)$ if $m=nd$ and $R_m(z)=0$ if $d\nmid m$, i.e., the sequence $\{R_m(z)\}$ is ``essentially" $\{P_n(z)\}$.
\end{rem}
Of all the results in the literature regarding the location of (the algebraic curve containing) the zeros of $\{P_n(z)\}$, the most general and crucial to us is {\cite[Theorem 1.1]{BO}} due to B\"ogvad, et. al. To state it, we need the following notation: For a non-empty subset $\mathcal{X}$ of nonzero complex polynomials, we set $\mathcal{Z}_{\mathcal{X}}$, (or $\mathcal{Z}_{A(z)}$ if $\mathcal{X}=\{A(z)\}$) to be the union of the set of zeros of all the polynomials in $\mathcal{X}$, i.e., $\displaystyle \mathcal{Z}_{\mathcal{X}}:=\{z^{\ast}\in \mathbb{C}: f(z^{\ast})=0 \text{ for some } f\in\mathcal{X}\}=\bigcup_{f(z)\in \mathcal{X}} \mathcal{Z}_{f(z)}$.
\begin{thm}[{\cite[Theorem 1.1]{BO}}]\label{conj:Tran} If $\{P_n(z)\}$ is the polynomial sequence generated by \eqref{uutiixxxx} and $z_0 \in \mathcal{Z}_{\{P_n(z)\}} - \mathcal{Z}_{A(z)}$, then
$z_0$ lies on the real algebraic curve $\displaystyle \Gamma:=\left\{z\in\mathbb{C}:{\rm Im}\left((-1)^k\frac{B^k(z)}{A^{\ell}(z)}\right)=0\right\}$.
\end{thm}
\begin{rem}\label{remark1.2}.
\begin{enumerate}
\item Theorem \ref{conj:Tran} generalizes the specific cases $k=2, 3$ and $4$ with $\ell =1$ proved in \cite[Theorems 1, 3, 5]{Tr1} respectively. 
\item The zeros of $P_n(z)$ become dense in $\displaystyle \mathcal{H}:=\Gamma\cap \Big\{z\in \mathbb{C}: 0\leqslant {\rm Re}\left((-1)^k\frac{B^k(z)}{A^{\ell}(z)}\right)< \frac{k^k}{(k-\ell)^{k-\ell}}\Big\}$ as $n$ tends to infinity for the cases $(k, \ell) \in \{(2,1),(3,1),(4,1)\}$, as proved in \cite[Theorems 1, 3, 5]{Tr1}. For $k\geqslant 5$ and $\ell=1$, Tran established the density of zeros of $P_n(z)$ in $\mathcal{H}$ for sufficiently large $n$, see \cite[Theorem 1]{Tr2} for details.
\item The case $\ell>1$ is of a different flavour, $($see Theorem \ref{in:conjectures}$)$ and will be established in Subsection \ref{B(ii)}.
\end{enumerate}
\end{rem}
Theorem \ref{conj:Tran} gives a partial characterisation of the real algebraic curve on which the zeros of the generated polynomials lie. This motivated the authors to search for the full characterisation of the curve $\Gamma$, as stated in Theorem \ref{in:conjectures}.
\begin{thm}\label{in:conjectures}
If $\{P_n(z)\}$ is the polynomial sequence generated by \eqref{uutiixxxx} with $\ell>1$ and $z_0 \in \mathcal{Z}_{\{P_n(z)\}} - \mathcal{Z}_{A(z)}$, then $z_0$ lies on the curve \begin{align*}
\Gamma_1:=\begin{cases}\left\{z\in\Gamma: {\rm Re}\left((-1)^{k-\ell}\frac{B^k(z)}{A^{\ell}(z)}\right)\leqslant 0\right\}, & \text{where }\mathcal{Z}_{\{B(z)\}} -\mathcal{Z}_{A(z)}\neq \emptyset\\ \left\{z\in\Gamma:{\rm Re}\left((-1)^{k-\ell}\frac{B^k(z)}{A^{\ell}(z)}\right) < 0\right\}, & \text{where }\mathcal{Z}_{\{B(z)\}} \subset \mathcal{Z}_{A(z)}\end{cases}
\end{align*}
\end{thm} 

Another problem that we deal with in the present paper concerns the location of zeros of trinomials. In general, the problem of the number and location of zeros of an arbitrary trinomial $at^k+bt^{\ell}+c$ (where $a,b,c\in\mathbb{C}-\{0\}$ and $k,\ell\in\mathbb{Z}$ such that $k\geqslant 3$ and $k>\ell\geqslant 1$) has a long history of study starting with J. Lambert in 1758, followed by L. Euler in 1777 and many others to the present, see \cite{AMelman, DBelkic} and the references therein. Of all the results in the literature about this problem, we are interested in those that take the form of bounds on the magnitude of the zeros or of sectors in the complex plane containing the zeros. These stem from J. Egerv\`ary's observation that, the zeros of trinomials can be interpreted as the equilibrium points of a force field created by unit masses that are located at the vertices of two regular concentric polygons centered at the origin in the complex plane. Utilising the symmetry and continuity properties of this force field, the roots of a trinomial can be separated according to their argument and moduli, see \cite{PGSzabo} for details. From these results, A. Melman in 2012 obtained finer results involving smaller annular sectors containing the zeros of a trinomial that take into account the magnitude of the coefficients. For further details, see \cite{AMelman}. 

Most of the approaches in the above studies could be characterised as being algebraic and/or geometric. However, most recently, in 2019, D. Belki\'c employed some analytic tools that primarily focus upon derivations of the analytical formulae for all the roots of trinomials through series developments using the Bell polynomials and the Fox-Wright function. These are based on the fact, that all the roots of the general $n$th degree trinomial admit certain convenient representations in terms of the Lambert and Euler series for the asymmetric and symmetric cases of the trinomial equation, respectively. As an application of their work, these analytical solutions are numerically illustrated in the genome multiplicity corrections for survival of synchronous cell populations after irradiation, see \cite{DBelkic} for details. 

Despite the above efforts, the problem of locating the zeros of a general trinomial is still unresolved. In the current paper, we take a different approach and instead study the location of ratios of zeros of trinomials. In particular, we choose an arbitrary $n\in \mathbb{Z}_{\geqslant 1}$, then a $z_0\in \mathcal{Z}_{P_n(z)}-\mathcal{Z}_{A(z)}$ and study the location of ratios of the zeros of the trinomil $D(t; z_0)=A(z_0)t^k+B(z_0)t^{\ell}+1$. We obtain the following results, namely: Theorem \ref{tn} and Corollary \ref{tn1}.

\begin{thm}[\textbf{Main result}]\label{tn}
If $\{P_n(z)\}$ is the polynomial sequence generated by \eqref{uutiixxxx} and $z_0 \in \mathcal{Z}_{\{P_n(z)\}} - \mathcal{Z}_{A(z)}$, then there exists at-least two zeros of $D(t; z_0)=A(z_0)t^k+B(z_0)t^{\ell}+1$ whose ratio is real and / or has modulus $1$.
\end{thm}
\begin{cor}\label{tn1}
If $\{P_n(z)\}$ is the polynomial sequence generated by \eqref{uutiixxxx} and $z_0 \in \mathcal{Z}_{\{P_n(z)\}} - \mathcal{Z}_{A(z)}$, then $D(t;z_0)$ has at-least two equimodular zeros and / or exactly three ``null-collinear"\footnote{Null-collinear means two points (for $k$ even) or three (for $k$ odd) points lie on a line via the origin $O$ with at-least one point on either side of $O$.} zeros for $k$ odd, or two null-collinear zeros for $k$ even.
\end{cor}

The remaining sections of the paper are devoted to proving Theorem \ref{in:conjectures} and Theorem \ref{tn} in reverse order. 
\section{Proofs}
\subsection{Some lemmata and important results} \label{B(i)}
Let $\mathcal{C}:=\{z\in \mathbb{C}:|z|=1\}$ and $h: X \to \mathbb{C}$, $\displaystyle w \mapsto \frac{(1-w^k)^k }{(1-w^\ell)^\ell(w^\ell-w^k)^{k-\ell}}$ where $X$ is the domain of $h$. 
\begin{lem}\label{go2}
If $z \in X\cap (\mathbb{R} \cup \mathcal{C})$, then $h(z) \in \mathbb{R}$. 
\end{lem}
\begin{proof}
The domain of $h$ is $X= \mathbb{C} -(\mu_\ell\cup\mu_{k-\ell}\cup\{0\})$ where $\mu_n$ denotes the set of $n$th complex zeros of unity. Now, there are two cases to consider, namely:
\begin{enumerate}[$(i)$]
\item if $z \in X \cap \mathbb{R}$, then $z\in \mathbb{R}$, hence $h(z) \in \mathbb{R}$ as $h$ is a quotient of two polynomials with real coefficients, i.e., $h\in \mathbb{R}(w)$.
\item if $z \in X \cap \mathcal{C}$, then  $h(z) \in \mathbb{R}$ by \cite[Lemma 2]{ndikus}. 
\end{enumerate}
\end{proof}
We now optimize $h$ over $X\cap \mathbb{R}$ since its optimal solutions on $\mathbb{R}$ have some connection with Theorem \ref{tn}. To do this, we first define the (analytic) extension $g$ of $h$ at $x=1$ as follows: $g:(X\cap\mathbb{R}) \cup \{1\}\to \mathbb{R}$,
\begin{align*} 
 x\mapsto \begin{cases} h(x), & x\in X\cap \mathbb{R},\\ \frac{k^k}{\ell^{\ell}(k-\ell)^{k-\ell}}, & x=1.\end{cases}
\end{align*}
To understand the behaviour of $g$, it suffices to study it over $[-1,1]-\{0\}$ since $h(q)=h(q^{-1})$. 
\begin{pro}
If
\begin{enumerate}[(i)]
\item $k$ is odd, then the global minimum value of $g$ is $\displaystyle g_{\rm min}=\frac{k^k }{\ell^\ell(k-\ell)^{k-\ell}}$.
\item $k$ is even, then the local minimum value of $g$ is $g_{\rm loc. min}=\displaystyle \frac{k^k }{\ell^\ell(k-\ell)^{k-\ell}}$ and local maximum value $g_{\rm loc. max}=0$.
\end{enumerate}
\end{pro}

\begin{proof}
The optimal value of $g$ occurs at the critical points of $g$, i.e., when either $g'(q)$ does not exist or $g'(q)=0$. Now
\begin{align}\label{KKK}
g'(q)=\frac{- (1-q^k)^{k-1}q^{\ell-1}(k-\ell + \ell q^k-k q^\ell)((k-\ell)q^{k} -k q^{k-\ell}+\ell)}{(1-q^\ell)^{\ell +1}(q^\ell- q^k)^{k-\ell +1}}. 
\end{align}
If $g'(q)$ does not exist, then $q\in (\mu_{\ell}\cup\mu_{k-\ell})\cap \mathbb{R}$. There are two cases to consider, depending on the parity of $k$:
\begin{enumerate}[(i)]
\item $k$ is odd: here $(\mu_{\ell}\cup\mu_{k-\ell})\cap \mathbb{R} = \{-1,1\}$, with $\displaystyle \lim_{q\to -1 }g(q)=\infty=\lim_{q\to 0 }g(q)$ and $\displaystyle g(1)=\frac{k^k }{\ell^\ell(k-\ell)^{k-\ell}}$. 
\item $k$ is even: here, $g'(q)$ does not exist if $q \in \{0,1\}$. Now $\displaystyle \lim_{q\to 0 }g(q)$ does not exist while $ \displaystyle g(1)=\frac{k^k }{\ell^\ell(k-\ell)^{k-\ell}}$.
\end{enumerate}
For $g'(q)=0$, we consider two cases depending on the parity of $k$ as follows:
\begin{enumerate}[(i)]
\item $k$ is odd: in this case, $g'(q)=0$ if and only if $q$ is a real solution to either 
\begin{align*}
k-\ell + \ell q^k-k q^\ell=0 \text{ or }(k-\ell)q^{k} -k q^{k-\ell}+\ell =0,
\end{align*}
as $q\notin \{0,1\}$. Since $k-\ell + \ell q^k-k q^\ell$ and $(k-\ell)q^{k} -k q^{k-\ell}+\ell $ are reciprocal polynomials of each other in the variable $q$, it suffices to solve any one of them. If $\gamma \neq 1$ is a solution to  $(k-\ell)q^{k} -k q^{k-\ell}+\ell =0$, then we have
\begin{align}\label{labe1}
(k-\ell){\gamma}^{k}=-\ell +k {\gamma}^{k-\ell}.
\end{align}
Evaluating $g$ at such a point $\gamma$ gives
\begin{align*}
g(\gamma) &= \frac{(1-\gamma^k)^k }{(1-\gamma^\ell)^\ell(\gamma^\ell-\gamma^k)^{k-\ell}}\\&=\frac{((k-\ell)-(k-\ell)\gamma^k)^k }{(1-\gamma^\ell)^\ell(1-\gamma^{k-\ell})^{k-\ell}\gamma^{\ell(k-\ell)}}\cdot \frac{1}{(k-\ell)^k}, ~~\mbox{}~ \hspace{2.6cm} by~~ \eqref{labe1}, 
\\&=\frac{(k-k\gamma^{k-\ell})^k }{(1-\gamma^\ell)^\ell(1-\gamma^{k-\ell})^{k-\ell}\gamma^{\ell(k-\ell)}}\cdot \frac{1}{(k-\ell)^k} \\&=\frac{(1-\gamma^{k-\ell})^\ell }{(1-\gamma^\ell)^\ell \gamma^{\ell(k-\ell)}}\cdot\frac{1}{(k-\ell)^\ell}\cdot \frac{k^k}{(k-\ell)^{k-\ell}},
\\&=\frac{(1-\gamma^{k-\ell})^\ell }{((k-\ell) \gamma^{k-\ell}-(k-\ell) \gamma^{k})^\ell}\cdot \frac{k^k}{(k-\ell)^{k-\ell}}, ~~\mbox{}~ \hspace{2.7cm} by~~\eqref{labe1}, 
\\&=\frac{(1-\gamma^{k-\ell})^\ell }{((k-\ell) \gamma^{k-\ell}+ \ell-k \gamma^{k-\ell})^\ell}\cdot\frac{k^k}{(k-\ell)^{k-\ell}}\\&=\frac{(1-\gamma^{k-\ell})^\ell }{(1-\gamma^{k-\ell})^\ell }\cdot\frac{k^k}{\ell^\ell(k-\ell)^{k-\ell}}=\frac{k^k}{\ell^\ell(k-\ell)^{k-\ell}}=g(1).
\end{align*}
Since $\displaystyle \lim_{q\to \pm \infty }g(q)=\infty $, and $g(q)>0$, it follows that $g(1)$ is a global minimum value of $g$ when $k$ is odd.

\item $k$ is even: in this case, $g'(q)=0$ if and only if either $q=-1$ or  $q\neq 1$ is a real solution to either 
\begin{align}\label{reciprocal2}
k-\ell + \ell q^k-k q^\ell=0 \text{ or }(k-\ell)q^{k} -k q^{k-\ell}+\ell =0.
\end{align}
Clearly, $g(-1)=0$ and for any $\gamma \neq 1$ which is a solution to either equation in \eqref{reciprocal2},
we have $g(\gamma)=g(1)$ by the calculation in $(i)$ above. Furthermore, since $g(q) \leqslant 0$ for $q <0$, we have $0$ as a local maximum. Similarly, as $g(q) \geqslant g(1) $ for $q >0$, we must have $g(1)$ as a local minimum of $g$.
\end{enumerate}
\end{proof} 
\begin{figure}[ht!]
\centering
\begin{subfigure}[b]{0.3\textwidth}
\includegraphics[height=4cm,width=5cm]{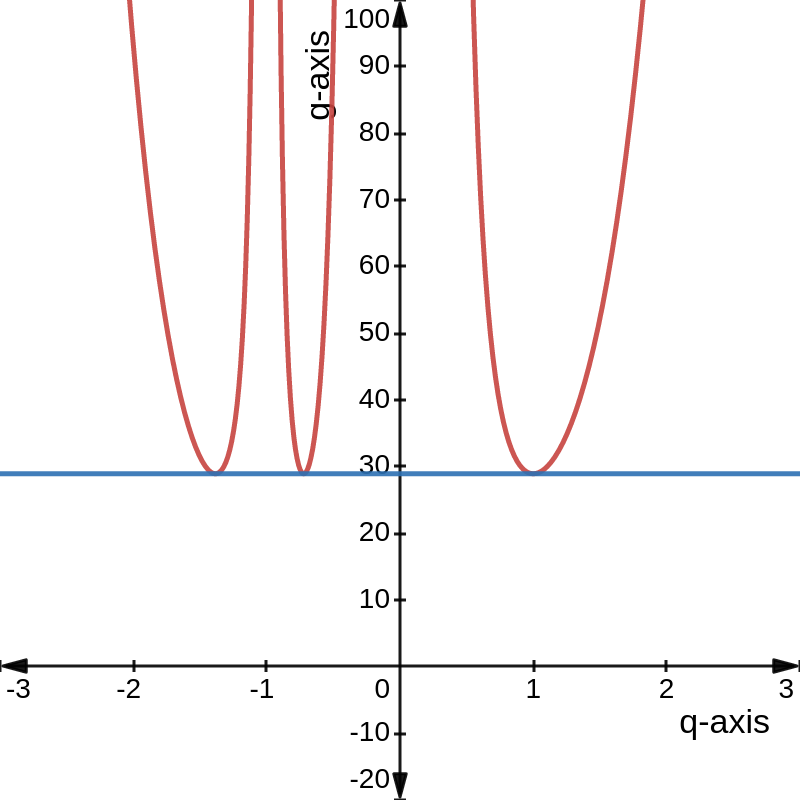}
\caption{Sketch for $g$ with $(k,\ell)=(5,3)$. \\ The blue line represents $g(1)$.}
\label{mandolo1}
\end{subfigure}
\hspace{2cm}
\begin{subfigure}[b]{0.3\textwidth}
\includegraphics[height=4cm,width=5cm]{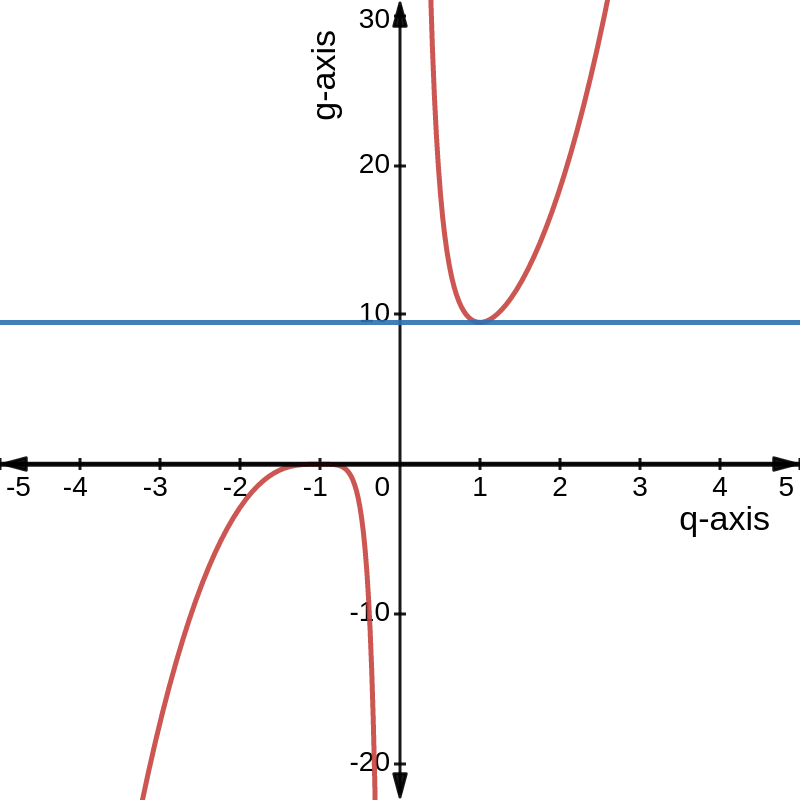}
\caption{Sketch for $g$ with $(k,\ell)=(4,1)$. \\ The blue line represents $g(1)$.}
\label{mandolo2}
\end{subfigure}
\caption{}
\label{man}
\end{figure}
\begin{defn}
Let $n\in \mathbb{Z}_{\geqslant 1}$. The $q$-discriminant of an $n$th degree polynomial $f(y)$ with zeros $y_1,\ldots,y_n$ and leading coefficient $a_n$ is; $($we use the convention that an empty product evaluates to 1$)$ \begin{align*}
\triangle_q(f(y))=a_n^{2n-2}q^{\frac{n(n-1)}{2}}\prod_{1\leqslant i<j\leqslant n}\left(y_i^2+y_j^2-\left(q+\frac{1}{q}\right)y_iy_j\right).
\end{align*}
\end{defn}
If $q=1$, then $\displaystyle \triangle_q(f(y))=a_n^{2n-2}\prod_{1\leqslant i<j\leqslant n}(y_i-y_j)^2=:\triangle(f(y))$, the classical discriminant of $f(y)$.
\begin{thm}[{\cite[Theorem 2]{ndikus}}]\label{ipt1} 
Let $k,\ell$ be coprime integers with $k>\ell \geqslant 1$ and $D(t;z)= A(z)t^k+ B(z)t^\ell+1$ where $A(z), B(z)\in \mathbb{C}[z]-\{0\}$. For any $z_0\in \mathbb{C}$ with $A(z_0)B(z_0)\neq 0$, the $q$-discriminant $\triangle_q(D(t;z_0))$ of the trinomial $D(t;z_0)$ is 
\begin{align} \label{In4ii}
(-1)^{\frac{k(k+3)}{2}}(A(z_0))^{k-1}(B(z_0))^{\ell-1}\frac{(1-q^k)^{k}}{(1-q)^{k}} \left(1-\frac{(1-q^\ell)^{\ell}(q^\ell -q^k)^{k-\ell}}{(q^k-1)^{k}}\frac{(B(z_0))^k}{(A(z_0))^\ell}\right).
\end{align}
\end{thm}
\begin{rem}\label{phd} 
We have the following remarks.
\begin{enumerate}
\item Vanishing of the $q$-discriminant for $D(t;z_0)$ at $q=1$ implies existence of repeated zeros for $D(t;z_0)$.
\item 
If $q_0\in \mathbb{C}$ is a ratio of zeros of $D(t;z_0)$, then $\triangle_{q_0}(D(t;z_0))=0$, i.e., $\displaystyle (-1)^k\frac{B^k(z_0)}{A^{\ell}(z_0)}=\frac{(1-q_0^k)^k }{(1-q_0^\ell)^\ell(q_0^\ell-q_0^k)^{k-\ell}}=:h(q_0)$. In addition, $\displaystyle h(q)=h(q_0)$ for any other ratio $q$ of distinct zeros of $D(t;z_0)$. 
\item If $z_{0}\in \mathcal{Z}_{\{P_n(z)\}} - \mathcal{Z}_{A(z)}$, then $\displaystyle (-1)^{k} \frac{B^k(z_0)}{A^\ell(z_0)}\in\mathbb{R}$ by \cite[Theorem 1.1]{BO}. 
\end{enumerate}
\end{rem}
\begin{rem}
Let $\{P_n(z)\}$ be the sequence of polynomials generated by \eqref{uutiixxxx}, $z_{0}\in \mathcal{Z}_{\{P_n(z)\}} - \mathcal{Z}_{A(z)}$, $ \displaystyle \alpha =(-1)^{k}\frac{B^k(z_0)}{A^\ell(z_0)}$, $\displaystyle \beta=\frac{k^k}{\ell^{\ell}(k-\ell)^{k-\ell}}$ and $\Omega_{k,\ell}(\alpha)$ be the number of ratios of distinct zeros of $D(t; z_{0})$ that are real. We have the following:
\begin{enumerate}[$(i)$]
\item If $k$ is odd, then 
\begin{align*}
\Omega_{k,\ell}(\alpha)=\begin{cases}
0, &\alpha<\beta,\\ 2, & \alpha=\beta,\\ 6, & \alpha>\beta.
\end{cases}
\end{align*}
In fact, for $\alpha=\beta$, there are three real ratios one of which is $q=1$ arising from non-distinct zeros of $D(t;z_0)$.
\item If $k$ is even, then 
\begin{align*}
\Omega_{k,\ell}(\alpha)=\begin{cases}
2, &\alpha <0, \alpha>\beta,\\ 1, & \alpha=0,\\ 0, & 0<\alpha\leqslant \beta.
\end{cases}
\end{align*}
In the case $\alpha=\beta$, there is actually a real ratio $q=1$ $($but it is not of distinct zeros of $D(t;z_0))$.
\end{enumerate}
\end{rem}
The $\beta$ is related to the separability threshold $\sigma(k,\ell)$ in \cite[Definition 2.3]{AMelman} as follows: $\beta=(\sigma(k,\ell))^k$. If $F(t):=t^k+at^{\ell}-1$ where $a\in\mathbb{C}-\{0\}$ as in \cite{AMelman}, then the reality of $\displaystyle (-1)^{k} \frac{B^k(z_0)}{A^\ell(z_0)}$ in Remark \ref{phd}(3), corresponds to the reality of $(-a)^k$. Corollary \ref{tn1} reveals a property exhibited by the roots of $F(t)$ if $(-a)^k$ is real, i.e., $F(t)$ has at-least two equimodular zeros and / or exactly three null-collinear zeros for $k$ odd, or exactly two null-collinear zeros for $k$ even. 

\begin{lem}\label{go21} 
If $\{P_n(z)\}$ is the sequence generated by \eqref{uutiixxxx}, $z_0 \in \mathcal{Z}_{\{P_n(z)\}} - \mathcal{Z}_{A(z)}$ and $q \in \mathbb{C}- (\mathbb{R}\cup \mathcal{C})$ is a ratio of distinct zeros of $D(t; z_0)$, then there is a ratio $q'\neq q$ of distinct zeros of $D(t; z_0)$ with $q' \in \mathbb{R}\cup \mathcal{C}$ and $\displaystyle h(q')=h(q)=(-1)^{k} \frac{B^k(z_0)}{A^\ell(z_0)}$.
\end{lem}
\begin{proof}
Suppose that $q \in \mathbb{C}- (\mathbb{R}\cup \mathcal{C})$ is a ratio of distinct zeros of $D(t; z_0)$. 
We show existence of another ratio $q'\in \mathbb{R}\cup \mathcal{C}$ of distinct zeros of $D(t;z_0)$. To this end, we let $A(z_0)=|A(z_0)|e^{\mathbf{i}\beta \pi}$ and $B(z_0)=|B(z_0)|e^{\mathbf{i}\theta \pi}$, where $\theta$ and $\beta$ are rational numbers. Since $k,\ell\in\mathbb{Z}_{\geqslant 1}$ are coprime, and $\displaystyle h(q)=(-1)^{k} \frac{B^k(z_0)}{A^\ell(z_0)}\in \mathbb{R}$, we have $\displaystyle \frac{k}{\ell}=\frac{\beta}{\theta} \pm \frac{\gamma}{\ell \theta}$ where $\gamma\in \{-1,0,1\}$. The transformation $Y=e^{\mathbf{i}\lambda \pi}t$ where $\displaystyle \lambda= \frac{\beta-\theta + \gamma}{k-\ell}$ transforms $D(t; z_0)$ into $\overline{D}(Y; z_0):= \pm |A(z_0)|Y^k\pm |B(z_0)|Y^\ell+1 \in \mathbb{R}[Y]$, (note, there are four possible polynomials). Via the  transformation $Y=e^{\mathbf{i}\lambda \pi}t$, the ratios of (distinct) zeros of the polynomials $D(t;z_0)$ and $\overline{D}(Y; z_0)$ are the same. So it is enough to study the ratios of distinct zeros of $\overline{D}(Y;z_0)$.
\begin{enumerate}[(i)]
\item When $k\geqslant 2$ is odd, then there must exist a real zero of $\overline{D}(Y;z_0)$. If it is the only zero, then there exists a complex conjugate pair $(y_1,\overline{y_1})$ of zeros of $\overline{D}(Y;z_0)$ whose ratio is in $\mathcal{C}$. In the event that there are more than one real zeros of $\overline{D}(Y;z_0)$, then we have a real ratio. 
\item When $k$ is even, either there are no real zeros or there must be at-least two real zeros of $\overline{D}(Y;z_0)$. In the latter case, there must exist a real ratio while in the former case, the zeros must occur in complex conjugate pairs, hence there must exist a ratio $q$ in $\mathcal{C}-\{1\}$.     
\end{enumerate}
\end{proof}
\begin{thm}\label{main}
If $\{P_n(z)\}$ is the polynomial sequence generated by \eqref{uutiixxxx}, $z_0 \in \mathcal{Z}_{\{P_n(z)\}} - \mathcal{Z}_{A(z)}$ and $q$ is a ratio of distinct zeros of $D(t; z_0)$, then $h(q)\in \mathbb{R}$ if and only if either $q \in (\mathbb{R}\cup\mathcal{C})\cap X$ or there is another ratio $q' \in (\mathbb{R}\cup\mathcal{C})\cap X$ of distinct zeros of $D(t; z_0)$.
\end{thm}

\begin{figure}[ht!]\centering
\includegraphics[scale=0.5]{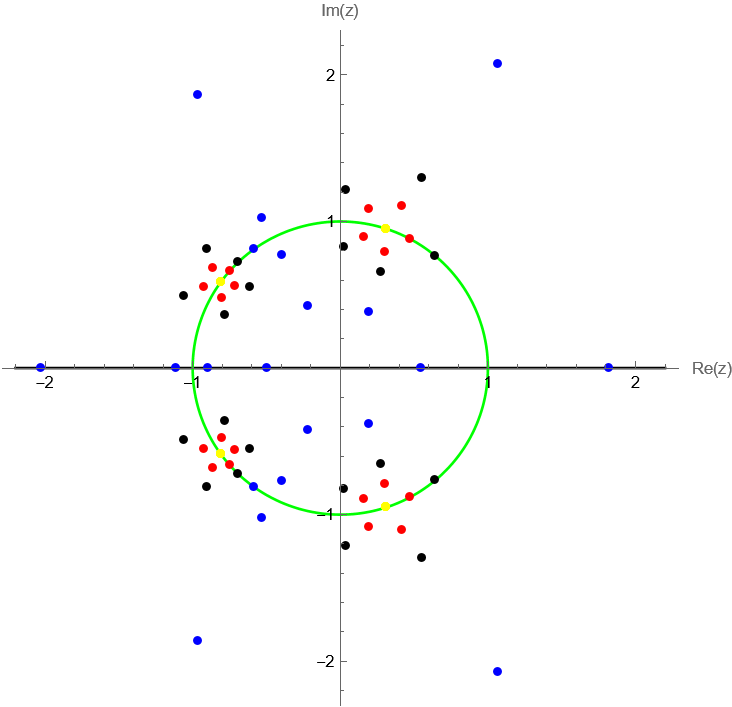}
\caption{The yellow, black, red, and blue dots correspond to the ratios of distinct zeros of $D(t; z_{17,0}), D(t; z_{23,0}), D(t; z_{56,0})$, and $D(t; z_{56,1})$ respectively. The green curve is $\mathcal{C}$, (unit circle centered at the origin).}
\label{a1}
\end{figure}

\begin{figure}[ht!]\centering
\includegraphics[scale=0.5]{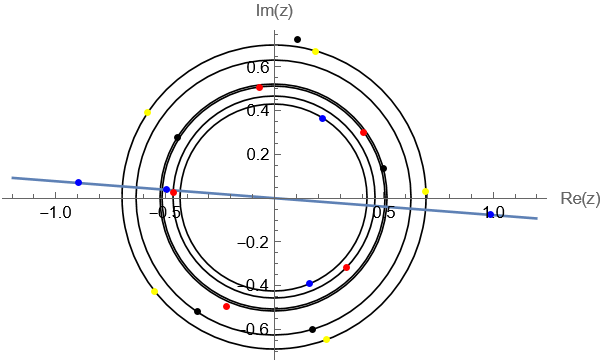}
\caption{The yellow, black, red, and blue dots correspond to the zeros of $D(t; z_{17,0}), D(t; z_{23,0}), D(t; z_{56,0})$, and $D(t; z_{56,1})$ respectively. The dots of the same color on the same circle represent equimodular zeros of a given $D(t;z_0)$. The three blue dots on a line via the origin are null-collinear.}
\label{a2}
\end{figure}

\begin{figure}[ht!]
\centering
\includegraphics[scale=0.5]{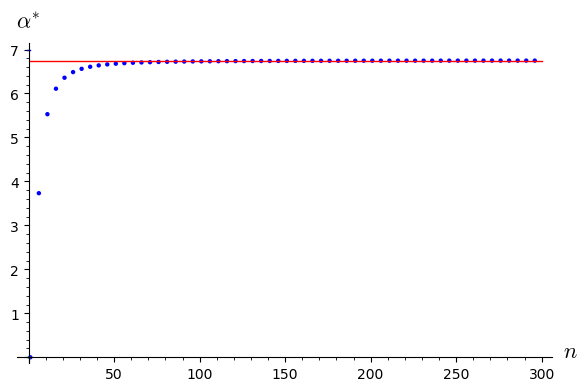}
\caption{Plot of $\displaystyle \alpha^{\ast}:=\max\left\{\frac{(-B(z_0))^3}{A(z_0)}:z_0\in \mathcal{Z}_{P_n(z)}-\mathcal{Z}_{A(z)}\right\}$ against $n\in \{1+5\lambda: \lambda=0,1,2,\ldots,59\}$ for $k=3$, $\ell=1$, $A(z)=2z$, and $B(z)=3z$.}
\label{limitingcase}
\end{figure}

\begin{proof}
We now show that if $q \in (\mathbb{R} \cup \mathcal{C})\cap X$ is a ratio of distinct zeros of $D(t; z_0)$, then $h(q)\in \mathbb{R}$. There are two cases: Either $q \in \mathbb{R}\cap X$, in which case $h(q) \in \mathbb{R}$ since $h(z)\in \mathbb{R}(z)$ or $q \in \mathcal{C}\cap X$, in which case we have $h(q) \in \mathbb{R}$ by Lemma \ref{go2}. On the other hand, if $q \notin X-(\mathbb{R} \cup \mathcal{C})\cap X$ is a ratio of distinct zeros of $D(t; z_0)$, then Lemma \ref{go21} implies existence of another ratio $q'$ of distinct zeros of $D(t; z_0)$ such that $q'\in X\cap(\mathbb{R}\cup \mathcal{C})$, which again implies that $h(q)=h(q')\in \mathbb{R}$.

$(\Rightarrow)$ Suppose that $h(q) \in \mathbb{R}$ where $q$ is a ratio of distinct zeros of $D(t; z_0)$. Since $q\neq 1$, Remark \ref{phd} implies that, \begin{align*}
h(q)=\frac{(1-q^k)^k }{(1-q^\ell)^\ell(q^\ell-q^k)^{k-\ell}}=(-1)^{k} \frac{B^k(z_0)}{A^\ell(z_0)}=\alpha.
\end{align*}
We consider the polynomial function 
\begin{align}\label{polyf}
f(w):=(1-w^k)^{k} -\alpha(1-w^\ell)^{\ell}(w^\ell -w^k)^{k-\ell}.
\end{align}
Since $\alpha\in \mathbb{R}[w]$, it follows that $f\in \mathbb{R}[w]$ and has degree $k^2$. Moreover, $\displaystyle f(w)= (1-w)^kH(w)$ where 
\begin{align*} H(w)= \left(\frac{1-w^k}{1-w}\right)^k-\alpha \left(\frac{1-w^\ell}{1-w}\right)^\ell\left(\frac{1-w^{k-\ell}}{1-w}\right)^{k-\ell}w^{\ell(k-\ell)}.
\end{align*}
It is clear that $H$ extends to $\widehat{H}\in \mathbb{R}[w]$ with $\displaystyle\widehat{H}(1):= \lim_{w \to 1} H(w)$. 

Clearly, $\displaystyle\widehat{H}(1)= \lim_{w \to 1} H(w)= k^k-\alpha \ell^\ell(k-\ell)^{k-\ell}=0$ if and only if $ \displaystyle \alpha= \frac{k^k}{\ell^\ell(k-\ell)^{k-\ell} }=g(1)$ which is forbidden (since $q\neq 1$, we consider distinct zeros of $D(t; z_0)$). Therefore, $1$ is a zero of $f$ of multiplicity $k$, hence $\deg(\widehat{H})= k^2-k\in 2 \mathbb{Z}_{\geqslant 1}$. Since the zeros of $\widehat{H}$ and $H$ are the same (up to multiplicity), without loss of generality, we can take $H=\widehat{H} \in \mathbb{R}[w]$. Now $\mathcal{W}_H\subset \mathcal{W}_f$ where $\mathcal{W}_H$ is the set of ratios of distinct zeros of $D(t; z_0)$ while $\mathcal{W}_f$ is the set of ratios of zeros of $D(t; z_0)$, not necessarily distinct. So if $q \in \mathcal{W}_H$, then so is $\overline{q}$ and ${q}^{-1}$ hence $\overline{{q}^{-1}}$ as $H(w) \in \mathbb{R}[w]$. Hence, $\displaystyle h(q)=h(\overline{q})=h(q^{-1})=h({\overline{q^{-1}}})=\alpha$.
Consequently, we have the following possibilities among the zeros of $H$.
\begin{enumerate}[(i)]
\item $q=\overline{q}$ implies that $ q\in \mathbb{R}- \{1\}$.
\item $q=q^{-1}$ implies that $q^2=1$ hence $q=-1 \in \mathbb{R}$ as $q\neq +1$.
\item  $q={\overline{q^{-1}}}$ implies that $q \overline{q}=1$ hence $|q|=1$, thence $q\in \mathcal{C} - (\mu_\ell\cup\mu_{k-\ell})$. 
\item Suppose that $q$ is a ratio of distinct zeros of $D(t; z_0)$, but none of the conditions (i)--(iii) holds. The conclusion that there is another ratio $q'$ of distinct zeros of $D(t; z_0)$ for which $q' \in (\mathbb{R} \cup \mathcal{C})\cap X$ follows from Lemma \ref{go21}.
\end{enumerate}
\end{proof}
Let us finally settle the main result of the present paper.
\begin{proof}[Proof of Theorem \ref{tn}] 
Let $z_0 \in \mathcal{Z}_{\{P_n(z)\}} - \mathcal{Z}_{A(z)}$, and $t_i:=t_i(z_0)$ for $i=1,2, \ldots, k$ be the zeros of 
\begin{align*}
D(t; z_0)= A(z_0)t^k+ B(z_0)t^\ell+1.
\end{align*}
There are two cases we shall consider; namely, repeated zeros and distinct zeros of $D(t; z_0)$.
\begin{enumerate}
\item \textit{Case 1:} (Repeated zeros). The classical discriminant of $D(t; z_0)$ vanishes. In this case, at least two of the zeros of $D$ are equal and nonzero. Therefore, at least one of the ratios of zeros of $D$ is $1 \in \mathbb{R}$ and on $\mathcal{C}$.
\item \textit{Case 2:} (Distinct zeros). In this case, the  classical discriminant of $D(t; z_0)$ does not vanish but its $q$-discriminant vanishes at any ratio $q$ of distinct zeros of $D(t; z_0)$. Let $q_0$ be a ratio of distinct zeros of $D(t; z_0)$, i.e., ($q_0 \neq 1$). Now, $\displaystyle h(q_0)=\frac{(1-q_0^k)^k }{(1-q_0^\ell)^\ell(q_0^\ell-q_0^k)^{k-\ell}}=(-1)^{k}\frac{B^k(z_0)}{A^\ell(z_0)}\in \mathbb{R}$ by \cite[Theorem 1.1]{BO}. Since $h(q_0)\in\mathbb{R}$, we have either $q_0\in (\mathbb{R}\cup\mathcal{C})\cap X$ or there is another ratio $q'\in (\mathbb{R}\cup\mathcal{C})\cap X$ of distinct zeros of $D(t;z_0)$, by Theorem \ref{main}.
\end{enumerate}
\end{proof} 
For a demonstration for Theorem \ref{tn}, we considered two arbitrary non-zero complex polynomials $A(z), B(z)$, with $\deg(A(z)B(z))\geqslant 1$ and coprime integers $k, \ell$ with $k >\ell \geqslant 1$. In particular, we considered (at random) the polynomials $A(z)=\mathbf{i}z^3 + z + 3\mathbf{i}, B(z)=z^2 - 2\mathbf{i} z + 7$, with $k=5, \ell=3,$ so that $k^2-k=20$. We numerically computed some few terms of the corresponding sequence $\{P_n(z)\}$, in particular $\{P_{17}(z),  P_{23}(z), P_{56}(z)\}$, and their corresponding zeros. Out of these zeros, we selected a few denoted by $z_{n,j}$ satisfying $A(z_{n,j})\neq 0$, computed $D(t;z_{n,j})$ as well as ratios of distinct zeros of $D(t;z_{n,j})$ and summarised the results in Table \ref{Bam2}. In Table \ref{Bam2}, we have 
\begin{align*} 
D(t; z_{17,0}):&=((25+4 \mathbf{i})-(22+2 \mathbf{i}) \sqrt{2}) t^5+1,\\D(t; z_{23,0}):&=-(8.85784 - 9.47542 \mathbf{i})  t^5 -(1.89645 - 3.91559 \mathbf{i}) t^3+1,\\D(t; z_{56,0}):&=(29.8505 + 14.9508 \mathbf{i}) t^5+ (3.50521 - 1.27839 \mathbf{i}) t^3+1,\\D(t; z_{56,1}):&=-(11.6461 + 4.7345 \mathbf{i}) t^5+ (11.1473 + 2.62973 \mathbf{i}) t^3+1.
\end{align*}
All computations were done using Mathematica \cite{Wolfram} and/or SageMath software \cite{SAGE} on a 64-bit machine with 5-cores and intel microprocessor with 4.0 GB of RAM and processor speed 2.4GHz.
\begin{table}[ht!]
\centering
\resizebox{13cm}{10cm}{
\begin{tabular}{|c|c|c|c|c|c|}
\hline
$n$ & $z_{n,j}\in \mathcal{Z}_{\{P_n(z)\}} - \mathcal{Z}_{A(z)}$ & $A(z_{n,j})=0$ & is $\displaystyle (-1)^{k}\frac{B^k(z_{n,j})}{A^\ell(z_{n,j})}\in \mathbb{R}?$ & $D(t; z_{n,j})$ & $q^\ast_{n,s}=q_{n,s} \in \mathbb{C}-\{0,1\}.$ \\ \hline
17 & $(1-2\sqrt{2})\mathbf{i}$ & no & yes, 0 & $D(t; z_{17,0})$ & ${\color{cyan}{-0.8090 \pm 0.5878\mathbf{i}}}$\\
 &  &  &  &  & ${\color{cyan}{0.3090 \pm 0.9511\mathbf{i}}}$\\ \hline
23 & $-0.6109-2.2046\mathbf{i}$ & no & yes, 0.71429  & $D(t; z_{23,0})$ & $0.5553 \pm 1.2966 \mathbf{i}$ \\
  &  &  &  &  & $0.0372 \pm 1.2114 \mathbf{i}$ \\ 
  &  &  &  &  & $-0.9024 \pm 0.8090\mathbf{i}$ \\ 
    &  &  &  &  & ${\color{cyan}{0.6444 \pm 0.7647\mathbf{i}}}$ \\ 
 &  &  &  &  & $0.0253 \pm 0.8247 \mathbf{i}$ \\
&  &  &  &  & $0.2791 \pm 0.6517 \mathbf{i}$ \\
 &  &  &  &  & ${\color{cyan}{-0.6900 \pm 0.7238\mathbf{i}}}$ \\
  &  &  &  &  & $-1.0553 \pm 0.4908\mathbf{i}$ \\
 &  &  &  &  & $-0.7791 \pm 0.3623 \mathbf{i}$ \\ 
&  &  &  &  & $-0.6144 \pm 0.5508 \mathbf{i} $\\ \hline
56 & $-0.2985+3.1410\mathbf{i}$ & no & yes, 0.01943 & $D(t; z_{56,0})$ & $0.4198 \pm 1.1044  \mathbf{i}$ \\
 &  &  &  &  & $0.1927 \pm 1.0854 \mathbf{i}$ \\ 
 &  &  &  &  & $0.1586 \pm 0.8932\mathbf{i}$ \\
&  &  &  &  & $0.3007 \pm 0.7911\mathbf{i}$ \\ 
  &  &  &  &   & ${\color{cyan}{0.4731 \pm 0.8810\mathbf{i}}}$ \\ 
 &  &  &  &  & $-0.8651 \pm 0.6833\mathbf{i}$ \\
&  &  &  &  & $-0.9198 \pm 0.5501 \mathbf{i}$ \\ 
 &  &  &  &  & ${\color{cyan}{-0.7475 \pm 0.6643\mathbf{i}}}$ \\
  &  &  &  &  & $-0.8007 \pm 0.4789 \mathbf{i}$\\ 
 &  &  &  &  & $-0.7119 \pm 0.5623 \mathbf{i}$ \\  \hline
56 &$ 1.9038+1.6907\mathbf{i}$ & no & yes, 99.18923 & $D(t; z_{56,1})$ & $1.0694 \pm 2.0772  \mathbf{i}$ \\
&  &  &  &  & $1.8190 $ \\ 
  &  &  &  &  & $0.5497 $ \\ 
 &  &  &  &  & $0.1960 \pm 0.3806\mathbf{i}$ \\
&  &  &  &  & $-0.9605 \pm 1.8655 \mathbf{i}$ \\ 
 &  &  &  &  & $-0.5280 \pm 1.0256 \mathbf{i}$ \\ 
&  &  &  &  & $-1.1134$ \\
 &  &  &  &  & ${\color{cyan}{-0.5810 \pm 0.8139\mathbf{i}}}$ \\ 
&  &  &  &  & $-0.8981 $\\
  &  &  &  &  & $-0.4937 $\\
 &  &  &  &  & $-0.2182 \pm 0.4237 \mathbf{i}$\\
&  &  &  &  & $-0.3968 \pm 0.7708 \mathbf{i}$ \\ &  &  &  &  & $-2.0254$ \\ \hline
\end{tabular}
}
\caption{Data with $k=5, \ell=3, A(z)=\mathbf{i}z^3 + z + 3\mathbf{i}, B(z)=z^2 - 2\mathbf{i} z + 7$.}\label{Bam2}
\end{table}
In Figure \ref{a1}, we illustrate Theorem \ref{tn} by showing the location of ratios of distinct zeros of the polynomials $D(t; z_{n,j})$ corresponding to the choice of zero $z_{n,j}$ of each of the polynomials $P_{17}(z), P_{23}(z)$ and $P_{56}(z)$ in Table \ref{Bam2}. In Figure \ref{a2}, we illustrate Corollary \ref{tn1}.
\begin{rem}
We have some remarks.
\begin{enumerate}
\item For each of the four colour points, i.e., blue, black, yellow and red in Figure \ref{a1}, there are at least two colour points on $\mathcal{C}$ which correspond to the highlighted points $($in pale blue$)$ in Table \ref{Bam2}.
\item For each non-real complex number $q^\ast_{n,j}$ lying off $\mathcal{C}$, Theorem \ref{tn} guarantees  existence of another ratio $q^{{\ast}'}_{n,j}$ of distinct zeros of $D(t; z_{n,j})$ such that $q^\ast_{n,j} \neq q^{{\ast}'}_{n,j}$ and either $q^{{\ast}'}_{n,j}$ is real or lies on $\mathcal{C}$.
\item The zero $z_{17,0}=(1-2 \sqrt{2})\mathbf{i}$ is a special zero for $P_{17}(z)$ since $B(z_{17,0})=0.$ If $B(z_0)=0$, then Theorem \ref{tn} trivially holds since the zeros of $\begin{displaystyle} D(t; z_0)=A(z_0)t^k+1\end{displaystyle}$ $($recall $A(z_0)\neq 0)$ are equimodular (of the form $\displaystyle t_j=\left|\frac{1}{A(z_0)}\right|e^{\frac{2\pi \mathbf{i} j}{k}}$ for $j=0,1,\ldots, k-1$). Clearly, any ratio of two distinct zeros of $D(t; z_0)$ lie on $\mathcal{C}$, and each occurs with multiplicity $k$. This explains why in Figure \ref{a1}, all the yellow colour points lie on $\mathcal{C}$ and there are four of them instead of 20. 
\end{enumerate}
\end{rem}
The ratios of distinct zeros of $D(t;z_0)$ may all be real only when $k^2-k\leqslant 6$, i.e., $k\in \{2,3\}$. For example, for $k=2$, we can consider $A(z)=z^2$, $B(z)=z^2-2z-5$ and $n=1$ which yields $\alpha=0$,  thence $q=-1$ as the only $($real$)$ ratio of distinct zeros. For $k=3$, a simple algebraic manipulation shows that $D(t;z_0)$, must have a pair of repeated zeros, implying $\displaystyle \alpha=\frac{27}{4}$ and hence exactly two real ratios of distinct roots of $D(t;z_0)$. Numerical experiments suggest that this might be a limiting case, with for example $A(z)=2z$, $B(z)=3z$, see Figure \ref{limitingcase} for the results.
\subsection{The real algebraic curve $\Gamma_1$.}\label{B(ii)}
In this subsection, we shall show that the inequalities stated in Theorem \ref{in:conjectures} are related to the zeros of $\{P_n(z)\}$. For completeness of the exposition, we consider the sequence $\{P_n(z)\}$ defined by the recurrence in Equation \eqref{uutiixxxx}. 
For the triple $(\ell,k,n)$, where $\ell,k, n\in \mathbb{Z}_{\geqslant 1}$ ($k$ and $\ell$ are fixed coprime integers with $k>\ell\geqslant 1$ used in the definition of the recurrence \eqref{uutiixxxx}), we define $S:=S_{\ell,k,n}=\{(i,j)\in\mathbb{Z}_{\geqslant 0}^2:i\ell+jk=n\}$. Clearly, $S$ is a finite set say with $s$ elements. This follows from the fact that, the linear Diophantine equation $x\ell+ky=n$ is a straight line with negative gradient, and the sought for solutions lie in the first quadrant of the $xy$-plane. In general, a linear Diophantine equation $x\ell+ky=n$ (with ${\rm gcd}(k,\ell)=1$) has integer solutions of the form $(x,y)=(x_0+ku,y_0-\ell u)$ where $x_0,y_0,u\in \mathbb{Z}$. Moreover, we can choose an $x_0$ and $y_0$ so that we parameterise $S$ as follows: $S=\{(i_u,j_u):=(x_0+ku,y_0-\ell u)\in\mathbb{Z}_{\geqslant 0}^2:u=1,\ldots,s\}$. 

To the triple $(\ell,k,n)$, we associate the polynomial $\displaystyle G_{\ell,k,n}(\tau):=\sum_{u=1}^{s}\binom{i_u+j_u}{i_u}\tau^{u-1}$ in $\tau$. This is the generating function for the number of north-east lattice paths starting from the origin to a lattice point $(i_u,j_u)$ lying on the linear Diophantine equation $x\ell+ky=n$ in the first quadrant of the $xy$-plane. Such polynomials $G_{\ell,k,n}(\tau)$ are real rooted since their coefficients form a P\'olya frequency sequence. In addition, all the roots of $G_{\ell,k,n}(\tau)$ are negative (as all its coefficients are non-negative), see \cite[Conjecture 1]{yamingyu} and \cite{BO} for details. In \cite[Lemma 2]{BO}, B\"ogvad et al., showed that 
\begin{align*}
P_{n}(z) = \pm B^{i_1}(z)A^{j_1}(z)G_{\ell,k,n}\left((-1)^{(k-\ell)}\frac{B^k(z)}{A^\ell(z)}\right).
\end{align*}
Using the real-rootedness property of $G_{\ell,k,n}(\tau)$, B\"ogvad et al., proved that, if $z_0 \in \mathcal{Z}_{P_n(z)}-\mathcal{Z}_{A(z)B(z)}$, then \begin{align*} G_{\ell,k,n}\left((-1)^{(k-\ell)}\frac{B^k(z_0)}{A^\ell(z_0)}\right)=0,
\end{align*}
hence $\displaystyle {\rm Im}\left((-1)^{(k-\ell)}\frac{B^k(z_0)}{A^\ell(z_0)}\right)=0$, as in Theorem \ref{conj:Tran}.

We pointed out in Remark \ref{remark1.2} that, if $\ell=1$, then $\displaystyle {\rm Re}\left((-1)^{(k-\ell)}\frac{B^k(z_0)}{A^\ell(z_0)}\right)$ is bounded in absolute value by $\displaystyle \frac{k^k}{(k-1)^{k-1}}$, but the case for $\ell>1$ had remained out of reach. We now address this case by showing that, if $\ell>1$, then $\displaystyle {\rm Re}\left((-1)^{(k-\ell)}\frac{B^k(z_0)}{A^\ell(z_0)}\right)$ can be arbitrarily large. To obtain this result using brute force, one would have to consider all the cases of $n$ and observe the behaviour of roots of $G_{\ell,k,n}(\tau)$. However, this is cumbersome and yields complicated expressions for one to deal with. A clever way is to use simpler cases in which the linear Diophantine equation has a solution on the axes, i.e., those in which either $n\equiv 0(\bmod\;\ell)$ or $n\equiv 0(\bmod\;k)$ and the simplest as $n\equiv 0(\bmod\;k\ell)$. For the remaining part of the arguments/computations, we consider zeros of the polynomials $P_n(z)$ where $n\equiv 0(\bmod\;k\ell)$. 
\begin{lem}\label{go21x1} 
If $n =\gamma k\ell$ for some $\gamma\in \mathbb{Z}_{\geqslant 1}$, then $G_{\ell,k,n}(\tau)$ is monic with constant term $1$ and sum of roots $\displaystyle -\binom{\frac{n}{\ell}-k+ \ell}{\frac{n}{\ell}-k}$.
\begin{proof}
If $n \equiv 0(\bmod\;k\ell)$, then among the non-negative solutions of the linear Diophantine equations $x\ell+ky=n$ are the pairs $\displaystyle \left(0, \frac{n}{k}\right)$ and $\displaystyle \left(\frac{n}{\ell}, 0\right)$. Since $S=\{(i_u,j_u):=(x_0+ku,y_0-\ell u)\in\mathbb{Z}_{\geqslant 0}^2:u=1,\ldots,s\}$, the first solutions corresponds to $u=1$ while the last to $u=s$. In particular, if $u=s$, then $j_s=0$ and $\displaystyle \binom{i_s+j_s}{i_s}=\binom{i_s}{i_s}=1$ which implies that $G_{\ell,k,n}(\tau)$ is monic. Similarly, if $u=1$, then $i_1=0$ and $\displaystyle \binom{i_1+j_1}{i_1}=\binom{j_1}{0}=1$ implying that $G_{\ell,k,n}(\tau)$ has constant term $1$. Since $G_{\ell,k,n}(\tau)$ is monic, the sum of its roots is $\displaystyle -\binom{i_{s-1}+j_{s-1}}{i_{s-1}}$, the negative of the coefficient of $\tau^{s-2}$ in $G_{\ell,k,n}(\tau)$. However, $\displaystyle \binom{i_{s-1}+j_{s-1}}{i_{s-1}}=\binom{\frac{n}{\ell}-k+\ell}{\frac{n}{\ell}-k}$ which follows from the choice of $u=s-1$ giving $j_{s-1}=\ell$ and $\displaystyle i_{s-1}=\frac{n}{\ell}-k$.
\end{proof}
\end{lem}
We determine $s$ explicitly. Let $N(a,b,n)$ be the number of non-negative integer solutions of a linear Diophantine equation $ax+by=c$ with ${\rm gcd}(a,b)=1$. In \cite[Theorem]{tripathi}, Tripathi showed that \begin{align*}  N(a,b,n)=\frac{n+a'(n)a+b'(n)b}{ab}-1,
\end{align*} where $a'(n)\equiv -na^{-1}(\bmod\;b)$ with $1\leqslant a'(n)\leqslant b$ and $b'(n)\equiv -nb^{-1}(\bmod\;a)$, with $1\leqslant b'(n)\leqslant a$.
\begin{lem}\label{lemmacountn}
If ${\rm gcd}(k,\ell)=1$ and $n=\gamma k\ell$ for some $\gamma\in \mathbb{Z}_{\geqslant 1}$, then $\displaystyle N(\ell,k,n)=\gamma+1$.
\end{lem}
\begin{proof}
Suppose $a=\ell$, $b=k$. Since ${\rm gcd}(k,\ell)=1$ and $n=\gamma k\ell$ for some $\gamma\in \mathbb{Z}_{\geqslant 1}$, it follows from \cite[Theorem]{tripathi} that $1\leqslant a'(n)\equiv 0(\bmod\;k)\leqslant k$, hence $a'(n)=k$ and $1\leqslant b'(n)\equiv 0(\bmod\;\ell)\leqslant \ell$, hence $b'(n)=\ell$. Therefore, 
\begin{align*} 
N(\ell,k,n)=\frac{\gamma k\ell+a'(n)a+b'(n)b}{\ell k}-1=\frac{\gamma k\ell+k\ell+\ell k}{k\ell}-1=\gamma+1.
\end{align*}
\end{proof}
\begin{lem}\label{ofono}
If $\ell>1$, ${\rm gcd}(k,\ell)=1$, $z_0 \in \mathcal{Z}_{P_n(z)}-\mathcal{Z}_{A(z)}$ and $\mathcal{Z}_{B(z)}-\mathcal{Z}_{A(z)}\neq \emptyset$, then $\displaystyle -\infty < {\rm Re}\left((-1)^{(k-\ell)}\frac{B^k(z_0)}{A^\ell(z_0)}\right) \leqslant 0$.
\end{lem}
\begin{proof}
Firstly, from the recurrence \eqref{uutiixxxx}, we observe that, if $n=\ell$, then $P_{\ell}(z)=-B(z)$. Therefore, if $z_0\in \mathcal{Z}_{P_{\ell}(z)}-\mathcal{Z}_{A(z)}$, then $B(z_0)=0$ hence $\displaystyle {\rm Re}\left((-1)^{(k-\ell)}\frac{B^k(z_0)}{A^\ell(z_0)}\right)=0$. Let $\tau_1,\ldots,\tau_{s-1}$ be the roots of $G_{\ell,k,n}(\tau)$ with $\tau^{\ast}=\min\{\tau_i\}_{i=1}^{s-1}$, i.e., $|\tau^{\ast}|$ is the largest root in absolute value. Since $G_{\ell,k,n}(\tau)$ is monic, it follows that $\displaystyle \sum_{i=1}^{s-1}\tau_i=- \binom{i_{s-1}+j_{s-1}}{i_{s-1}}$. Since all the roots of $G_{\ell,k,n}(\tau)$ are negative, $\displaystyle \frac{1}{s-1} \sum_{i=1}^{s-1}\tau_i\geqslant \frac{1}{s-1} \sum_{i=1}^{s-1}\tau^{\ast}=\tau^{\ast}$ hence 
\begin{align*}|\tau^{\ast}|=-\tau^{\ast}\geqslant -\frac{1}{s-1} \sum_{i=1}^{s-1}\tau_i=\frac{1}{s-1}\binom{i_{s-1}+j_{s-1}}{i_{s-1}}.
\end{align*}
If $n=\gamma k\ell$, for some $\gamma\in \mathbb{Z}_{\geqslant 1}$, then $s=\gamma+1$ by Lemma \ref{lemmacountn}, hence \begin{align*} 
\frac{1}{s-1}\binom{i_{s-1}+j_{s-1}}{i_{s-1}}=\frac{1}{\gamma}\binom{\frac{n}{\ell}-k+\ell}{\frac{n}{\ell}-k}=\frac{1}{\gamma}\binom{(\gamma-1)k+\ell}{(\gamma-1)k}\approx \frac{(\gamma-1)^{\ell}k^{\ell}}{\ell!\gamma}.
\end{align*}As $\gamma$ hence $n$ tends to $\infty$, it follows that $\displaystyle \frac{1}{s-1} \binom{i_{s-1}+j_{s-1}}{i_{s-1}}$, (hence $|\tau^{\ast}|$) becomes arbitrarily large.
\end{proof}

\begin{proof}[Proof of Theorem \ref{in:conjectures}]
If $\ell>1$, and $\mathcal{Z}_{B(z)}\subset \mathcal{Z}_{A(z)}$, then for any $z_0\in \mathcal{Z}_{P_n(z)}-\mathcal{Z}_{A(z)}$, we have $B(z_0)\neq 0$, hence $\displaystyle {\rm Re}\left((-1)^{(k-\ell)}\frac{B^k(z_0)}{A^\ell(z_0)}\right) \neq 0$. Therefore, the complete proof Theorem \ref{in:conjectures} follows from Theorem \ref{conj:Tran} + Lemma \ref{ofono}.
\end{proof}
\begin{rem}
In the situation $(k, \ell)=(k,1)$, we have 
\begin{align*}
\frac{1}{s-1}\binom{i_{s-1}+j_{s-1}}{i_{s-1}}=\frac{1}{\gamma}\binom{(\gamma-1)k+1}{(\gamma-1)k}= \frac{(\gamma-1)k+1}{\gamma}\end{align*} 
which implies that the average values of the roots of $G_{k,\ell,n}(\tau)$ are bounded by $k$ for all $n\equiv 0(\bmod\:k\ell)$. This is in agreement with the results obtained earlier in \cite[Theorems 1, 3, 5]{Tr1} respectively and \cite[Theorem 1]{Tr2} where it was proved that any $z_0 \in \mathcal{Z}_{P_n(z)}-\mathcal{Z}_{A(z)}$ must satisfy the condition that $\displaystyle 0\leqslant {\rm Re}\left((-1)^k\frac{B^k(z_0)}{A(z_0)}\right)\leqslant \frac{k^k}{(k-1)^{k-1}}$. This is because, as $n$ tends to $\infty$ so does $\gamma$ and hence the expression $\displaystyle \frac{(\gamma-1)k+1}{\gamma}$ tends to $\displaystyle k < \frac{k^k}{(k-1)^{k-1}}$ for any $k\in \mathbb{Z}_{\geqslant 2}$.
\end{rem}

\section*{\refname}
\nocite{*}
\bibliographystyle{elsarticle-num-names}
\bibliography{bibratiosofdistinctzerosnew}
\appendix
\section{SageMath code}
\begin{lstlisting}[basicstyle=\ttfamily\scriptsize, language=Python]
#Sage code for generating the data in Table 1, and graphs in Figure 1 and Figure 2.
#Given coprime integers k and l such that k>l>=1, arbitrary complex polynomials A and B such that degree of A(z)B(z) is atleast 1 and positive integer n, the code computes the polynomial Pn(z) in polynomial sequence {Pn(z)} generated by the recurrence in equation (1). Furthermore, it computes the zeros of Pn(z). Using a zero z0 of one's choice (such that A(z0) is nonzero), the code numerically computes the trinomial D(t;z0) = A(z0)t^k+B(z0)t^l+1, its zeros, ratios of zeros of D(t;z0) and then plots both the zeros and ratios of zeros of D(t;z0) on two separate graphs. 
#To minimise the numerical errors (mostly when computing the ratios of zeros of D(t;z0)), we include a precision parameter, and for comparison, we have included another way of computing these ratios using alpha = (-1)^k(B(z0))^k/(A(z0))^l and polynomial f, equation (6) in the proof of Theorem 2.8. 

#variables
z,t,q=var('z,t,q') #names of variables of different polynomials used in the code.

#parameters (needed to reproduce the graphs in the paper)
k, l = 5, 3  # coprime integers k and l such that k>l>=1
A, B = i*z^3+z+3*i, z^2 -2*i*z+7  # complex polynomials A and B such that degree of A(z)B(z) is >= 1
precision = 50 # precision paramater, change to 100, 150 or 200 for accuracy (at the expense of speed)
c=[(1,1,0),(0,0,0),(1,0,0),(0,0,1)] # colors (yellow, black, red and blue) for different polynomials.

#data used in the code to reproduce the graphs, includes choice of n and label for zero zj of Pn(z).
list1=[(17,1),(23,2),(56, 9), (56,35)] # (n,j)=(23,2) corresponds to P23(z) and z2 = -0.6109 - 2.2046i
#one must choose a label j and therefore a zero zj such that A(zj) is not zero.

#functions
def recursive(n,k,l):  # inputs parameters n,k and l as described above and computes Pn(z) recursively.
    if 1-k<=n<0:  
        return 0  # returns polynomials in initial conditions.
    elif n==0:
        return 1 # returns polynomial P0(z).
    else: 
        return -B*recursive(n-l,k,l)-A*recursive(n-k,k,l) # returns polynomial Pn(z) for n>=1.
def wzeros(w, precision):  # returns complex zeros of a function w(z) computed numerically.
    if w.degree(z)>=1:
        wroots=w.roots(ring=ComplexField(precision)) 
        return wroots  # complex zeros of w(z) computed to given precision.
    elif w==0: # for this choice of w, some parts of the code may not work well.
        return w, "is a zero polynomial." 
    else:  # for this choice of w, some parts of the code may not work well.
        return w, "is a nonzero constant polynomial, therefore, it has no zeros."  
def polyfforratios(q,x): # polynomial in the proof of Theorem 2.8 (its zeros will be used as controls)
    return (1-q^k)^k-x*(1-q^l)^l*(q^l-q^k)^(k-l)

#our plots (initialisation with reference unit circles colored green)
s1=plot(circle((0,0),1,aspect_ratio=1,rgbcolor=(0,1,0),linestyle='-'))
s11=plot(circle((0,0),1,aspect_ratio=1,rgbcolor=(0,1,0),linestyle='-'))
count=0 # counter as we run through list1, will help in plotting all graphs on the same plot.

for j in list1: #program to run through list1 to generate the graphs
    n, choiceofzero=j[0], j[1] # e.g., (23,2)
    w = (recursive(n,k,l)).expand() # expands the recursive polynomial as a polynomial in z
    print("Our choice of Pn(z) is: ", w) # e.g., P23(z)
    zerosofw = wzeros(w,precision) # e.g., this list contains zeroes of P23(z)
    print("The zeros of our choice of Pn(z) with their multiplicities are: ")
    for counter, s in enumerate(zerosofw): # e.g., label, zero of P23(z), multiplicities
        print(" %d: "% counter, zerosofw[counter][0], zerosofw[counter][1])
    yourchoice = zerosofw[choiceofzero][0] # e.g., z2=-0.6109 - 2.2046i, for P23(z)
    w0 = A(z=yourchoice)*t^k+B(z=yourchoice)*t^l+1 # we already checked that A(zj) is nonzero,
    # e.g., D(t;z2)=-(8.85784 - 9.47542i)t^5 - (1.89645 - 3.91559i)t^3 + 1.
    print("Your choice of zero is: ", yourchoice, "and the corresponding D(t,z0) is: ", w0)
    zerosofw0 = w0.roots(ring=ComplexField(precision)) # seros of D(t;z0)
    print("The zeros of D(t,z0) with their multiplicities are:")
    for counter, s in enumerate(zerosofw0):  # e.g., label, zero of D(t,z2), multiplicities
        print(" %d: "% counter, zerosofw0[counter][0], zerosofw0[counter][1])
    alpha = (-1)^k*(B(z=yourchoice))^k/(A(z=yourchoice))^l # alpha parameter to be used in f
    print(alpha) # e.g., alpha=
    w1 = polyfforratios(q,alpha) # polynomial f with alpha, see (6) in the proof of Theorem 2.8
    ratiosofzerosofD = w1.roots(ring=ComplexField(precision))  
    print("The theoretical ratios of (not necessarily distinct) zeros of D(t,z0):")
    for counter, s in enumerate(ratiosofzerosofD): # e.g., label, zero of f, multiplicities
        print(" %d: "% counter, ratiosofzerosofD[counter][0], ratiosofzerosofD[counter][1])
    zerosofDnum = [b[1][0] for b in enumerate(zerosofw0)] # e.g., zeros of D(t,z2) computed numerically
    cart_zerosofDnum = cartesian_product([zerosofDnum,zerosofDnum]) # e.g., cartesian product.
    ratiosofzerosofDnum = [b[0]/b[1] for b in cart_zerosofDnum] # e.g., ratios of zeros of D(t,z2)
    print("The numerical ratios of (not necessarily distinct) zeros of D(t,z0):")
    for counter, s in enumerate(ratiosofzerosofDnum):# label, ratios of zeros of D(t,z0), multiplicities
        print(" %d: "% counter, ratiosofzerosofDnum[counter])
    points_ratiosofzerosofD=[d[0] for d in ratiosofzerosofD] # ratios of zeros of D(t,z0) to be plotted 
    s1=s1+points(points_ratiosofzerosofD, rgbcolor=c[count]) # append plot of ratios of zeros of D(t,z0)
    points_zerosofD=[d[0] for d in zerosofw0] # zeros of D(t,z0) to be plotted below
    s11=s11+points(points_zerosofD, rgbcolor=c[count]) # append plot of ratios of zeros of D(t,z0)
    for j in points_zerosofD: #appending circles centered at origin containing each zero of D(t,z0)
        s11=s11+plot(circle((0,0),j.abs(),aspect_ratio=1,thickness=0.3,rgbcolor=(0,0,0),linestyle='-'))
    count=count+1
s1.axes_labels([r"$\mathrm{Re}(z)$","$\mathrm{Im}(z)$"])
save(s1,'/tmp/figure1.png') 
s11.axes_labels([r"$\mathrm{Re}(z)$","$\mathrm{Im}(z)$"])
save(s11,'/tmp/figure2.png')
os.system('display /tmp/figure2.png')
os.system('display /tmp/figure2.png')
show(s1) # for jupyter notebook, sagemath cell or cocalc interpreters.
show(s11) # for jupyter notebook, sagemath cell or cocalc interpreters.
\end{lstlisting}
\end{document}